\documentclass[12pt]{article}

\usepackage{amsmath, amssymb}
\usepackage{amsfonts, mathrsfs}
\usepackage[cp1251]{inputenc}
\usepackage[active]{srcltx}
\usepackage[ukrainian, russian, english]{babel}

\usepackage{amsfonts, amssymb, mathrsfs}

\begin{document}

\selectlanguage{ukrainian}
\thispagestyle{empty}

 \noindent {\small  {{{{УДК 517.5 }}}}     } \vskip 1.5mm

 \noindent \textbf{А.\,С. Сердюк, І.\,В. Соколенко} {\small (Iн-т математики НАН України, Київ)} \vskip 2mm
 \noindent{serdyuk@imath.kiev.ua, sokol@imath.kiev.ua}
 \vskip 2mm
\noindent  \textbf{Наближення інтерполяційними тригонометричними поліномами\\ в метриках просторів \boldmath{$L_p$} на класах періодичних цілих функцій
}
 \vskip 3.5mm

\noindent  {\small Встановлено асимптотичні рівності для точних верхніх меж наближень інтерполяційними тригонометричними поліномами з рівномірними розподілом вузлів інтерполяції $x_k^{(n-1)}=\frac{2k\pi}{2n-1},\ k\in\mathbb{Z},$ у метриках просторів $L_p, 1\le p\le\infty,$ на класах $2\pi$-періодичних функцій, які зображуються у вигляді згорток функцій $\varphi, \ \varphi\perp1,$ що належать одиничній кулі з простору $L_1$ із фіксованими твірними ядрами, у яких модулі коефіцієнтів Фур'є $\psi(k)$ задовольняють умову $\lim\limits_{k\rightarrow\infty} \psi(k+1)/\psi(k)=0.$ Аналогічні оцінки встановлено і на класах $r$-диференційовних функцій $W^r_1$ при швидко зростаючих показниках гладкості $r$ $(r/n\rightarrow\infty)$.}

 \vskip 3.5mm
 \noindent \textbf{А.\,С. Сердюк, И.\,В. Соколенко}  \vskip 2mm
 
 \noindent  \textbf{Приближение интерполяционными тригонометрическими полиномами\\ в метриках пространств \boldmath{$L_p$} на классах периодических целых функций
 }
 \vskip 3.5mm
 
\noindent   {\small Получены асимптотические равенства для точных верхних граней приближений  интерполяционными тригонометрическими полиномами с равномерным распределением узлов интерполяции $x_k^{(n-1)}=\frac{2k\pi}{2n-1},\ k\in\mathbb{Z},$ в метриках пространств $L_p, 1\le p\le\infty,$ на классах $2\pi$-периодических функций, представимых в виде сверток функций $\varphi, \ \varphi\perp1,$ принадлежащих единичному шару пространства $L_1$, с фиксированными производящими ядрами, у которых модули коэффициентов Фурье  $\psi(k)$ удовлетворяют условию $\lim\limits_{k\rightarrow\infty} \psi(k+1)/\psi(k)=0.$ Аналогичные оценки получены и на классах $r$-дифференцируемых функций $W^r_1$ при быстро возрастающих показателях гладкости $r$ $(r/n\rightarrow\infty)$.}
 
 \vskip 3.5mm
 \noindent \textbf{A.\,S. Serdyuk, I.\,V. Sokolenko}  \vskip 2mm
 
 \noindent  \textbf{Approximation by interpolation trigonometric polynomials\\ in metrics of the spaces \boldmath{$L_p$} on the classes of periodic entire functions
 }
 \vskip 3.5mm

\noindent   {\small We obtain the asymptotic equalities for the least upper bounds
	of approximations  by interpolation trigonometric polynomials with the equidistant nodes $x_k^{(n-1)}=\frac{2k\pi}{2n-1},\ k\in\mathbb{Z},$  in metrics of
	the spaces $L_p$ on classes of $2\pi$-periodic functions, representable as convolutions of functions $\varphi, \ \varphi\perp1,$ which belongs to the unit ball of the space $L_1$, and fixed generating kernels in the case where modules of their Fourier coefficients $\psi(k)$ satisfy the condition $\lim\limits_{k\rightarrow\infty} \psi(k+1)/\psi(k)=0.$ We obtain similar estimates on the classes of $r$-differentiable functions $W^r_1$ for the quickly increasing exponents of smoothness $r$ $(r/n\rightarrow\infty)$.
	}

 \newpage
 \setcounter{page}{1}
 \normalsize

Нехай  $\ L_p,\ 1\le p<\infty,$ --- простір $2\pi$-періодичних сумовних
у $p$-му степені на $[0, 2\pi)$  функцій $f$ зі стандартною нормою $
\|f\|_{p}=\bigg(\int\limits_{0}^{2\pi}|f(t)|^pdt\bigg)^{1/p};
$ $\ L_{\infty}$ --- простір $2\pi$-періодичних вимірних і суттєво
обмежених функцій $f$ з нормою $
\|f\|_{\infty}=\mathop {\rm ess \sup}\limits_{t}
|f(t)|;$ $C$ --- простір $2\pi$-періодичних неперервних функцій
$f$, у якому норма задається рівністю $\|f\|_{C}=\max\limits_{t}|f(t)|.
$

Позначимо через $C^\psi_{\bar{\beta},1}$ класи $2\pi$-періодичних функцій $f$, які зображуються у вигляді згортки
\begin{equation}\label{1}
f(x)=\frac{a_0}2+\frac1\pi\int\limits_0^{2\pi}\Psi_{\bar{\beta}}(x-t)\varphi(t)dt, \ \ a_0\in\mathbb{R}, 
\end{equation}
в якій $\varphi\perp1,\  \|\varphi\|_1\le1,$ а $\Psi_{\bar{\beta}}(\cdot)$~--- ядра вигляду
\begin{equation}\label{2}
\Psi_{\bar{\beta}}(t)=\sum_{k=1}^{\infty}\psi(k)\cos\left(kt-\frac{\pi\beta_k}{2}\right),\ \ \ \beta_k\in\mathbb{R},\ \ \ \psi(k)>0,
\end{equation}
\begin{equation}\label{3}
\sum_{k=1}^{\infty}\psi(k)<\infty.
\end{equation}
 
 Якщо послідовності $\bar{\beta}=\{\beta_k\}_{k=1}^\infty$ є стаціонарними послідовностями, тобто $\beta_k=\beta,\ k\in\mathbb{N},$ $\ \beta\in\mathbb{R},$  то класи $C^\psi_{\bar{\beta},1}$ позначатимемо через $C^\psi_{{\beta},1}$.

Якщо $\psi(k)=k^{-r}, r>1,$ то класи $C^\psi_{\bar{\beta},1}$ та $C^\psi_{{\beta},1}$ позначатимемо, відповідно, через $W^r_{\bar{\beta},1}$ і $W^r_{{\beta},1}$. Останні класи є відомими класами Вейля-Надя. Якщо $r\in\mathbb{N}$ i $\beta=r$, то класи $W^r_{{\beta},1}$ є класами $W^r_1$~--- $2\pi$-періодичних функцій, що мають абсолютно неперервні похідні до $(r-1)$-го порядку включно і такі, що їх $r$-та похідна належить одиничній кулі простору $L_1$ (тобто $\|f^{(r)}\|_1\le1).$

Якщо  $\psi(k)=e^{-\alpha k^r}, \ \alpha>0,\ r>0,$ класи $C^\psi_{\bar{\beta},1}$ і $C^\psi_{{\beta},1}$ будемо позначати через $C^{\alpha,r}_{\bar{\beta},1}$  та $C^{\alpha,r}_{{\beta},1}$, відповідно. Останні класи називають іноді класами узагальнених інтегралів Пуассона.

У даній роботі розглядатимуться класи $C^\psi_{\bar{\beta},1}$ за умови, що послідовність $\psi(k)>0$ така, що 
\begin{equation}\label{4}
\lim\limits_{k\rightarrow\infty} \frac{\psi(k+1)}{\psi(k)}=0.
\end{equation}
Як випливає з [\ref{Stepanets2002_1}, c.~139-145], класи $C^\psi_{\bar{\beta},1}$ за виконання умови (\ref{4}) складаються із функцій, які допускають регулярне продовження в усю комплексну площину, тобто складаються із цілих функцій. З іншого боку, як показано в [\ref{StepanetsSerdyukShydlich2008}, c.~1703], для того, щоб функція $f$  належала до множини усіх дійснозначних на дійсній осі цілих функцій, необхідно і достатньо, щоб вона могла бути зображена згорткою вигляду (\ref{1}) у якій $\varphi\in L_1,$ а коефіцієнти $\psi(k)$ ядра $\Psi_{\bar{\beta}}$ вигляду (\ref{2})  задовольняли умову (\ref{4}).

Нехай $f\in C.$ Через $\ \tilde S_{n-1}(f; x)\ $  позначатимемо
тригонометричний поліном порядку $\ n-1,\ $ що інтерполює $\ f(x)\ $ у  рівномірно розподілених вузлах $\ x_k^{(n-1)}=\frac{2k\pi}{2n-1},\ $ $k\in \mathbb{Z},$ тобто такий, що 
$$
\tilde S_{n-1}(f; x_k^{(n-1)}) = f(x_k^{(n-1)}), \ \ \ k\in \mathbb{Z}.
$$

Порядкові оцінки збіжності інтерполяційних поліномів $\tilde{S}_{n-1}(f; \cdot)$ до $f$ в метриках просторів $C$ i $L_p$, що виражались в термінах послідовностей найкращих наближень функцій в $C$ i $L_p$, одержані у роботах [\ref{Sharapudinov1983}, \ref{Oskolkov1986}].

Розглянемо величину
\begin{equation}\label{5}
 \tilde{\cal E}_n(C^\psi_{\bar{\beta},1})_{L_p}=\sup\limits_{f\in C^\psi_{\bar{\beta},1}}\|f(\cdot)-\tilde{S}_{n-1}(f; \cdot)\|_p.
\end{equation}

У даній роботі будуть встановлені асимптотично точні оцінки величин (\ref{5}) при $n\rightarrow\infty$ для довільних $1\le p\le\infty,\ \beta_k\in\mathbb{R}\ $ i $\psi(k)$ таких, що задовольняють умову (\ref{4}).

При $p=1$ асимптотична поведінка величин вигляду (\ref{5}) при $n\rightarrow\infty$ в залежності від тих чи інших обмежень на послідовності $\psi(k)$ та $\beta_k$ досліджувалась у роботах [\ref{Motornyj1990}--\ref{Serdyuk2002}]. Зокрема у [\ref{Serdyuk2000}, c.~994] за виконання умови (\ref{4}) для довільних $\beta_k\in\mathbb{R}$ встановлено асимптотичну рівність 
\begin{equation}\label{6}
\tilde{\cal E}_n(C^\psi_{\bar{\beta},1})_{L_1}=\frac{16}{\pi^2}\psi(n)+O(1)\left(\frac{\psi(n)}{n}+\sum_{k=n+1}^{\infty}\psi(k)\right),
\end{equation}
в якій $O(1)$ рівномірно обмежена відносно усіх розглядуваних параметрів.

Крім того, у роботі [\ref{SerdyukVojtovych2010}, c.~279-280] отримано результати, з яких випливає, що за виконання умови (\ref{4}) при довільних $\beta_k\in\mathbb{R}$ має місце асимптотична рівність
\begin{equation}\label{7}
\tilde{\cal E}_n(C^\psi_{\bar{\beta},1})_{L_\infty}=\frac2{\pi}\psi(n)+O(1)\sum_{k=n+1}^{\infty}\psi(k),
\end{equation}
в якій $O(1)$ --- величина, рівномірно обмежена відносно усіх розглядуваних параметрів.

Питання про асимптотичну поведінку величин $\tilde{\cal E}_n(C^\psi_{\bar{\beta},1})_{L_p}$, $\beta_k\in\mathbb{R}$, за виконання умови (\ref{4}) при $1<p<\infty$ залишалось відкритим. 

Основним результатом є наступне твердження.

\textbf{Теорема 1.}\textit{  Нехай $1\le p\le\infty,$  $\ \beta_k
\in\mathbb{R},\ $  а  $\ \psi(k)>0$ задовольняє умову $(\ref{3}).$
Тоді при всіх $n\in\mathbb{N}$ має місце оцінка
\begin{equation}\label{8}
\tilde{\cal E}_n(C^\psi_{\bar{\beta},1})_{L_p}=\frac{2^{1-\frac1p}}{\pi^{1+\frac1p}}\|\cos t\|^2_p\psi(n)+O(1)\left(\frac{\psi(n)}{n}+\sum\limits_{\nu=n+1}^\infty\psi(\nu)\right),
\end{equation}
в якій $O(1)$ --- величина, рівномірно обмежена відносно усіх розглядуваних параметрів. Якщо, крім того, $\psi(k)$ задовольняє умову $(\ref{4})$, то оцінка $(\ref{8})$ є асимптотичною при $n\rightarrow\infty$ рівністю.}

\textbf{ \textit{Доведення.} } Згідно з формулою (9) роботи [\ref{Serdyuk2000}] для довільної функції $f$ з класу $C^\psi_{\bar{\beta},1}$ в кожній точці $x\in\mathbb{R}$ виконується рівність
\begin{equation}\label{9}
\tilde\rho_n(f;x)=f(x)-\tilde{S}_{n-1}(f;x)=
$$
$$
=\frac2\pi\psi(n)\sin\frac{2n-1}2x\int\limits_{0}^{2\pi}\sin\left(nt-\frac x2+\frac{\pi\beta_n}2\right)\varphi(t)dt+\tilde{r}_{n+1}(f;x),
\end{equation}
в якій 
\begin{equation}\label{10}
\tilde{r}_{n+1}(f;x):=\frac1\pi\int\limits_{0}^{2\pi}\sum_{\nu=n+1}^{\infty}\psi(\nu)\left(\cos\left(\nu(t-x)+\frac{\pi\beta_\nu}{2}\right)-\bar{\omega}_\nu(x;t)\right)\varphi(t)dt,
\end{equation}
де функції $\bar{\omega}_\nu(x;t), \nu=n,n+1,\ldots$
 означаються за допомогою формул
\begin{equation}\label{11}
\bar{\omega}_{m(2n-1)+k}(x;t)=\cos\left(m(2n-1)t+k(t-x)+\frac{\pi\beta_{m(2n-1)+k}}{2}\right),
$$
$$
 m\in\mathbb{N},\ \ k=0,\pm1,\ldots,\pm(n-1).
\end{equation} 
Із (\ref{10}) i (\ref{11}) випливає, що для довільних $1\le p\le\infty$ i $f\in C^\psi_{\bar{\beta},1}$ має місце оцінка
\begin{equation}\label{12}
\|\tilde{r}_{n+1}(f;\cdot)\|_p\le\frac2\pi\left\|\int\limits_{0}^{2\pi}\sum_{\nu=n+1}^{\infty}\psi(\nu)|\varphi(t)|dt\right\|_p\le
\frac{2^{1+\frac1p}}{\pi^{1-\frac1p}}\sum_{\nu=n+1}^{\infty}\psi(\nu)\le 2\sum_{\nu=n+1}^{\infty}\psi(\nu).
\end{equation}
Згідно з формулами (\ref{9}) i (\ref{12}) для довільних $1\le p\le\infty, \bar{\beta}=\{\beta_k\}_{k=1}^\infty, \beta_k\in\mathbb{R},$ i $\psi(k)>0,$ що задовольняють умову (\ref{3}), виконується оцінка
\begin{equation}\label{13}
\tilde{\cal E}_n(C^\psi_{\bar{\beta},1})_{L_p}=\frac2\pi\psi(n)A_n(p)+O(1)\sum_{\nu=n+1}^{\infty}\psi(\nu),
\end{equation}  
в якій 
\begin{equation}\label{14}
A_n(p)=A_n(p;\bar{\beta}):=\sup\limits_{
	\begin{array}{c}
	\|\varphi\|_1\le1\\
	\varphi\perp1
	\end{array}}
\left\|\sin \frac{2n-1}2x\int\limits_{0}^{2\pi}\sin\left(nt-\frac x2+\frac{\pi\beta_n}2\right)\varphi(t)dt\right\|_p.
\end{equation}
Дослідимо асимптотичну поведінку величин $A_n(p), \ 1\le p\le\infty,$ при $n\rightarrow\infty.$

З цією метою встановимо істинність наступної двосторонньої оцінки:
\begin{equation}\label{15}
\frac12\|\Phi_{n,\pi\beta_n}(x)\|_p\le A_n(p)\le\frac12\sup_{\theta\in\mathbb{R}}\|\Phi_{n,\theta}(x)\|_p, \ n\in\mathbb{N},\ 1\le p\le\infty,
\end{equation}
в якій
\begin{equation}\label{16}
\Phi_{n,\theta}(x):=\cos\left(nx-\frac\theta2\right)g_\theta(x)+\sin\left(nx-\frac\theta2\right)h_\theta(x),\ \ \theta\in\mathbb{R},
\end{equation}
\begin{equation}\label{17}
g_\theta(x):=1-\cos(x-\theta),
\end{equation}
\begin{equation}\label{18}
h_\theta(x):=-\sin(x-\theta).
\end{equation}
Для знаходження необхідної оцінки зверху величини $A_n(p)$ скористаємось узагальненою нерівністю Мінковського
\begin{equation}\label{19}
\left\|\int\limits_0^{2\pi}f(\cdot,u)du\right\|_p\le\int\limits_0^{2\pi}\left\|f(\cdot,u)\right\|_pdu,\ \ \ 1\le p\le\infty,
\end{equation}
(див. [\ref{Kornejchuk1987}, c. 395]).

Згідно з (\ref{14}) i (\ref{19})
\begin{equation}\label{20}
A_n(p)\le\sup\limits_{
	\begin{array}{c}
	\|\varphi\|_1\le1\\
	\varphi\perp1
	\end{array}}
\int\limits_{0}^{2\pi}\left\|\sin \frac{2n-1}2x\sin\left(nt-\frac x2+\frac{\pi\beta_n}2\right)\right\|_p\varphi(t)dt\le
$$
$$
\le \sup_{\theta\in\mathbb{R}}\left\|\sin \frac{2n-1}2x\sin\left(-\frac x2+\frac{\theta}2\right)\right\|_p.
\end{equation}
Оскільки при будь-якому $\theta\in\mathbb{R}$
\begin{equation}\label{21}
\sin \frac{2n-1}2x\sin\left(-\frac x2+\frac{\theta}2\right)=\frac12\left(\cos\left(nx-\frac{\theta}2\right)-\cos\left((n-1)x+\frac{\theta}2\right)\right)=
$$
$$
=\frac12\left(\cos\left(nx-\frac{\theta}2\right)(1-\cos(x-\theta))-\sin\left(nx-\frac{\theta}2\right)\sin(x-\theta)\right)=\frac12\Phi_{n,\theta}(x),
\end{equation}
то в силу (\ref{20}) i (\ref{21})
\begin{equation}\label{22}
A_n(p)\le\frac12\sup_{\theta\in\mathbb{R}}\|\Phi_{n,\theta}(x)\|_p.
\end{equation}

Оцінимо знизу величину $A_n(p).$ Для цього розглянемо при кожному $n\in\mathbb{N}$ і достатньо малому $\delta>0\ (\delta<\frac\pi n)$ $2\pi$-періодичну функцію $\varphi_{n,\delta}(t)$ таку, що на $[-\frac{\delta}{2}, 2\pi-\frac{\delta}{2}]$ задається за допомогою рівностей
\begin{equation}\label{23}
\varphi_{n,\delta}(t)=
\left\{
\begin{array}{rl}
\frac1{2\delta}, & t\in(-\frac{\delta}{2знам},\frac{\delta}{2знам})\\
-\frac1{2\delta}, & t\in(\frac{\pi}{nзнам}-\frac{\delta}{2знам},\frac{\pi}{nзнам}+\frac{\delta}{2знам})\\
0, & t\in[-\frac{\delta}{2знам}, 2\pi-\frac{\delta}{2знам}]\setminus\{(-\frac{\delta}{2знам},\frac{\delta}{2знам})\cup(\frac{\pi}{nзнам}-\frac{\delta}{2знам},\frac{\pi}{nзнам}+\frac{\delta}{2знам})\}.\\
\end{array}
\right.
\end{equation}
Оскільки $\|\varphi_{n,\delta}\|_1\le1$ i $\varphi_{n,\delta}\perp1,$ то, як випливає з (\ref{14})
\begin{equation}\label{24}
A_n(p)\ge
\left\|\sin \frac{2n-1}2x\int\limits_{0}^{2\pi}\sin\left(nt-\frac x2+\frac{\pi\beta_n}2\right)\varphi_{n,\delta}(t)dt\right\|_p.
\end{equation}
В силу (\ref{23})
\begin{equation}\label{25}
\int\limits_{0}^{2\pi}\sin\left(nt-\frac x2+\frac{\pi\beta_n}2\right)\varphi_{n,\delta}(t)dt=\frac1{2\delta}\int\limits_{-\frac{\delta}{2знам}}^{\frac{\delta}{2знам}}\sin\left(nt-\frac x2+\frac{\pi\beta_n}2\right)dt-
$$
$$
-\frac1{2\delta}\int\limits_{\frac{\pi}{nзнам}-\frac{\delta}{2знам}}^{\frac{\pi}{nзнам}+\frac{\delta}{2знам}}\sin\left(nt-\frac x2+\frac{\pi\beta_n}2\right)dt=\frac{1}{n\delta}\left(\cos\left(-\frac x2+\frac{\pi\beta_n}2-\frac{n\delta}{2знам}\right)-\cos\left(-\frac x2+\frac{\pi\beta_n}2+\frac{n\delta}{2знам}\right)\right).
\end{equation}
Із  (\ref{24}) i (\ref{25}) одержуємо нерівність
\begin{equation}\label{26}
A_n(p)\ge
\left\|\sin \frac{2n-1}2x \ \frac{\cos\left(-\frac x2+\frac{\pi\beta_n}2-\frac{n\delta}{2знам}\right)-\cos\left(-\frac x2+\frac{\pi\beta_n}2+\frac{n\delta}{2знам}\right)}{n\delta}\right\|_p.
\end{equation}
Обираючи $\delta$ настільки малим, щоб $\delta=o(\frac1n)$, і переходячи до границі  в правій частині нерівності (\ref{26}) при $n\delta\rightarrow0,$ з урахуванням формули (\ref{21}), застосованої при $\theta=\pi\beta_n,$ одержимо
\begin{equation}\label{27}
A_n(p)\ge
\left\|\sin \frac{2n-1}2x \sin\left(-\frac x2+\frac{\pi\beta_n}2\right)\right\|_p=\frac12 \|\Phi_{n,\pi\beta_n}(x)\|_p.
\end{equation}
Оцінка (\ref{15}) є наслідком формул (\ref{22}) і (\ref{27}).

Для знаходження асимптотичної рівності для величин $\|\Phi_{n,\theta}(\cdot)\|_p$ при $n\rightarrow\infty$ і довільних фіксованих $\theta\in\mathbb{R}$ i $1\le p\le\infty$ скористаємось наступним твердженням роботи [\ref{Serdyuk2005}, c. 1083].

\textbf{Лема 1.}\textit{  Нехай $1\le p\le\infty$ і $2\pi$-періодичні функції $g(x)$ i $h(x)$ мають обмежену варіацію на $[0,2\pi]$, якщо $p=1,$ або належать класу Гельдера $KH^1=\{f\in C: |f(x+\delta)-f(x)|\le K\delta,\  x,\delta\in\mathbb{R} \},$ якщо $1<p\le\infty.$ Тоді для функції 
$$
\Phi(x)=g(x)\cos(nx+\alpha)+h(x)\sin(nx+\alpha),\ \ \alpha\in\mathbb{R},\ \ n\in\mathbb{N},
$$
справджуються асимптотичні формули
\begin{equation}\label{28}
\left.
\begin{array}{c}
\|\Phi\|_p\\
\inf\limits_{c\in\mathbb{R}}\|\Phi(\cdot)-c\|_p\\
\sup\limits_{\theta\in\mathbb{R}}\frac12\|\Phi(\cdot+\theta)-\Phi(\cdot)\|_p
\end{array}
\right\}=
\frac{\|\cos t\|_p}{(2\pi)^{1/p}знам}\|r\|_p+O(1)\frac Mn,
\end{equation} 
в яких
\begin{equation}\label{29}
r(t)=\sqrt{g^2(x)+h^2(x)},
\end{equation}
\begin{equation}\label{30}
M=M_p=\left\{ 
\begin{array}{ll}
\mathop{V}\limits_{-\pi}^\pi(g)+\mathop{V}\limits_{-\pi}^\pi(h), & \mbox{при}\ p=1,\\
K+p^{-1}\|r\|_p^{1-p}\mathop{V}\limits_{-\pi}^\pi(r^p),& \mbox{при}\ 1<p<\infty,\\
K,& \mbox{при}\ p=\infty,
\end{array}
\right.
\end{equation}
а величина $O(1)$ рівномірно обмежена відносно усіх розлядуваних параметрів.
}

Покладемо в термінах леми 1 $g(x)=g_\theta(x),\ h(x)=h_\theta(x),\ \alpha=-\frac\theta2,\ \Phi(x)=\Phi_{n,\theta}(x).$ Тоді згідно з (\ref{17}), (\ref{18}) i (\ref{29}) при будь-якому $\theta\in\mathbb{R}, 1\le p\le\infty,$
\begin{equation}\label{31}
r(x)=r_\theta(x)=\sqrt{g^2_\theta(x)+h^2_\theta(x)}=\sqrt{2(1-\cos(x-\theta))}=2\left|\sin\frac{x-\theta}{2знам}\right|,
\end{equation}
\begin{equation}\label{32}
\mathop{V}\limits_{-\pi}^\pi(r^p_\theta)=2^{p+1},
\end{equation}
\begin{equation}\label{33}
\|r_\theta\|_p=2\left\|\sin\frac{x-\theta}{2знам}\right\|_p=2\|\cos t\|_p.
\end{equation}
Крім того, згідно з (\ref{30})--(\ref{33}) при $p=1$
\begin{equation}\label{34}
M_p=M_1=\mathop{V}\limits_{-\pi}^\pi(g_\theta)+\mathop{V}\limits_{-\pi}^\pi(h_\theta)=8,
\end{equation}
при   $1<p<\infty$
\begin{equation}\label{35}
M_p=K+\frac{\mathop{V}\limits_{-\pi}^\pi(r^p_\theta)}{p\|r_\theta\|_p^{p-1}}=1+\frac{4}{p\|\cos t\|_p^{p-1}},
\end{equation}
при $p=\infty$
\begin{equation}\label{36}
M_p=M_\infty=K=1.
\end{equation}
Оскільки при довільних $1<p<\infty$
$$
p\|\cos t\|_p^p=p\|\sin t\|_p^p\ge4p\int\limits_{0}^{\pi/2}\left(\frac2\pi t\right)^pdt=2\pi\frac p{p+1}\ge\pi
$$
i
$$
\|\cos t\|_p\le2\pi,
$$
то з рівностей (\ref{35}) отримуємо оцінку
\begin{equation}\label{37}
M_p\le9,\ \ 1<p<\infty.
\end{equation}

Застосувавши лему 1 до функції $\Phi(x)=\Phi_{n,\theta}(x)$  і врахувавши формули (\ref{33}), (\ref{34}), (\ref{36}) i (\ref{37}), із (\ref{28}) отримуємо рівномірну відносно усіх параметрів оцінку
\begin{equation}\label{38}
\|\Phi_{n,\theta}(\cdot)\|_p=\frac{2^{1-\frac1p}}{\pi^{\frac1p}}\|\cos t\|_p^2+O(1)\frac1n,\ \ \ 1\le p\le\infty,\ \ \theta\in\mathbb{R}.
\end{equation}

Із формул (\ref{13}), (\ref{15}) i (\ref{38}) випливає, що за умови (\ref{3}) справджується оцінка (\ref{8}).

Залишається показати, що при виконанні рівності (\ref{4}) оцінка (\ref{8}) перетворюється 	в асимптотичну при $n\rightarrow\infty$ рівність.

Розглянемо величину
\begin{equation}\label{42}
\varepsilon_n=\varepsilon_n(\psi):=\sup\limits_{k\ge n}\frac{\psi(k+1)}{\psi(k)}.
\end{equation}
В силу (\ref{4}) і (\ref{42}) величина $\varepsilon_n$	 монотонно спадає до нуля. А, отже, при $n$ таких, що $\varepsilon_n<1$
\begin{equation}\label{43}
\sum\limits_{k=n+1}^\infty\psi(k)=\psi(n)\frac{\psi(n+1)}{\psi(n)}+\psi(n)\frac{\psi(n+1)}{\psi(n)}\frac{\psi(n+2)}{\psi(n+1)}+\psi(n)\frac{\psi(n+1)}{\psi(n)}\frac{\psi(n+2)}{\psi(n+1)}\frac{\psi(n+3)}{\psi(n+2)}+\ldots=
$$
$$
=\psi(n)\sum\limits_{k=0}^\infty\prod_{j=0}^{k}\frac{\psi(n+j+1)}{\psi(n+j)}\le\psi(n)\sum\limits_{k=0}^\infty\prod_{j=0}^{k}\varepsilon_{n+j}\le
\psi(n)\sum \limits_{k=0}^\infty\varepsilon_{n}^{k+1}=\psi(n)\frac{\varepsilon_n}{1-\varepsilon_n}.
\end{equation}
З урахуванням співвідношень (\ref{43}), при $n$ таких, що  	$\varepsilon_n<1,$ з (\ref{8}) випливає рівномірна по усіх параметрах оцінка
\begin{equation}\label{44}
\tilde{\cal E}_n(C^\psi_{\bar{\beta},1})_{L_p}=\psi(n)\left(\frac{2^{1-\frac1p}}{\pi^{1+\frac1p}}\|\cos t\|^2_p +O(1)\left(\frac{1}{n}+ \frac{\varepsilon_n}{1-\varepsilon_n}\right)\right),
\end{equation}
яка при $n\rightarrow\infty$ є асимптотичною рівністю.

Теорему доведено.

При $p=1$	 формула (\ref{8}) перетворюється у формулу (\ref{6}), а при $p=\infty$ формула (\ref{8}) випливає з (\ref{7}).

Наведемо наслідок з теореми 1 у випадку, коли $\psi(k)=e^{-\alpha k^r},\ \alpha>0,\ r>1.$

\textbf{Теорема 2.}\textit{  Нехай $r>1,\ \alpha>0,\ 1\le p\le\infty,$  $\ \beta_k
	\in\mathbb{R}.\ $	Тоді для всіх номерів $n$ таких, що 
\begin{equation}\label{46}
n^{1-r}\ln (n+1)\le\alpha r,
\end{equation}
 має місце рівномірна відносно усіх розглядуваних параметрів оцінка
\begin{equation}\label{47}
\tilde{\cal E}_n(C^{\alpha,r}_{\bar{\beta},1})_{L_p}=e^{-\alpha n^r}\left(\frac{2^{1-\frac1p}}{\pi^{1+\frac1p}}\|\cos t\|^2_p+O(1)\frac1{n}\right).
\end{equation}
}	

\textbf{ \textit{Доведення.} } Для функції $\psi(k)=e^{-\alpha k^r},\ \alpha>0,\ r>1,$ величина $\varepsilon_n$ вигляду (\ref{42}) оцінюється наступним чином:
\begin{equation}\label{48}
\varepsilon_n=e^{-\alpha( (n+1)^r-n^r)}\le e^{-\alpha r n^{r-1}}.
\end{equation}
В силу (\ref{46})
\begin{equation}\label{49}
e^{-\alpha r n^{r-1}}\le \frac1{n+1}.
\end{equation}
Із  (\ref{48}) і  (\ref{49}) випливає нерівність
\begin{equation}\label{50}
\frac{\varepsilon_n}{1-\varepsilon_n}\le \frac2n.
\end{equation}
Формула (\ref{47}) випливає із (\ref{44}) і (\ref{50}).

Застосуємо також оцінку (\ref{8}) теореми 1 до класів $W^r_{\bar{\beta},1},$ тобто у випадку, коли $\psi(k)=k^{-r},$ $r>1.$ Оскільки 
$$
\sum\limits_{k=n+1}^\infty\frac1{k^r}<\frac1{(n+1)^r}+\int\limits_{n+1}^{\infty}\frac{dt}{t^r}= \frac1{(n+1)^r}\left(1+\frac{n+1}{r-1}\right),
$$
то, враховуючи 	монотонне зростання до $e^{-1}$ послідовності $\left(1-\frac1{n+1}\right)^{n+1}$, одержуємо
\begin{equation}\label{51}
n^r\sum\limits_{k=n+1}^\infty\frac1{k^r}<\left(\frac n{n+1}\right)^r\left(1+\frac{n+1}{r-1}\right)=\left(1-\frac1{n+1}\right)^{(n+1)\frac r{n+1}}\left(1+\frac{n+1}{r-1}\right)<
$$
$$
<e^{-\frac r{n+1}}\left(1+\frac{n+1}{r-1}\right).
\end{equation}
Із (\ref{51}) випливає рівномірна по всіх параметрах оцінка
\begin{equation}\label{52}
n^r\sum\limits_{k=n+1}^\infty\frac1{k^r}=O(1)e^{-\frac r{n+1}}\left(1+\frac{n}{r-1}\right), \ \ \ r>1.
\end{equation}

Як наслідок, з теореми 1 отримуємо таке твердження.

\textbf{Теорема 3.}\textit{  Нехай $1\le p\le\infty,\ r>1,\ $  $\ \beta_k
	\in\mathbb{R},\ n\in\mathbb{N}.$	Тоді має місце рівномірна по всіх параметрах оцінка
	\begin{equation}\label{53}
	\tilde{\cal E}_n(W^r_{\bar{\beta},1})_{L_p}=\frac1{n^r}\left(\frac{2^{1-\frac1p}}{\pi^{1+\frac1p}}\|\cos t\|^2_p+O(1)\left(\frac1{n}+e^{-\frac r{n+1}}\left(1+\frac{n}{r-1}\right)\right)\right).
	\end{equation}
}

Зазначимо, що при $r\ge n+1$ оцінка (\ref{53}) набуває вигляду
	\begin{equation}\label{54}
\tilde{\cal E}_n(W^r_{\bar{\beta},1})_{L_p}=\frac1{n^r}\left(\frac{2^{1-\frac1p}}{\pi^{1+\frac1p}}\|\cos t\|^2_p+O(1)\left(\frac1{n}+e^{-\frac r{n+1}}\right)\right).
\end{equation}

У випадку, коли $\frac rn\rightarrow\infty$ при $n\rightarrow\infty$, формула (\ref{54}) є асимптотичною рівністю.

При $p=1$ i $\beta_k=\beta,\ \beta\in\mathbb{R},\ r\ge n+1$ рівність (\ref{54}) запишеться у вигляді
\begin{equation}\label{54'}
\tilde{\cal E}_n(W^r_{{\beta},1})_{L_1}=\frac1{n^r}\left(\frac{16}{\pi^2} +O(1)\left(\frac1{n}+e^{-\frac r{n+1}}\right)\right).
\end{equation}

Оцінка (\ref{54'}) є інтерполяційним аналогом результатів робіт С.Б. Стєчкіна [\ref{Stechkin1980}, теорема 4] та С.О. Теляковського [\ref{Telyakovskii1989}], в яких досліджувалась асимптотика величин
$$
{\cal E}_n(W^r_{{\beta},1})_{L_1}=\sup\limits_{f\in W^r_{{\beta},1}}\|f(\cdot)-S_{n-1}(f; \cdot)\|_p,
$$
де $S_{n-1}(f;\cdot)$ --- частинна сума Фур'є порядку $n-1$ функції $f$.

Формула (\ref{54'}) доповнює результат В.П. Моторного [\ref{Motornyj1990}], згідно з яким при довільних $r\in\mathbb{N}$ виконується нерівність
\begin{equation}\label{55}
\tilde{\cal E}_n(W^r_1)_{L_1}\le\frac{2K_{r-1}\ln n}{\pi n^r} +\frac{O(1)}{n^r},
\end{equation}
де $K_r$ --- константи Фавара:
$$
K_m=\frac4\pi\sum\limits_{\nu=0}^\infty\frac{(-1)^{\nu(m+1)}}{(2\nu+1)^{m+1}},\ \ \ m=0,1,\ldots
$$
При цьому у випадку $r=2$ у формулі (\ref{55}) можна поставити знак "дорівнює"\,, тобто справджується асимптотична рівність
$$
\tilde{\cal E}_n(W^2_1)_{L_1}=\frac1{n^2}(\ln n+O(1)).
$$

\footnotesize
\begin{enumerate}

\item\label{Stepanets2002_1} {\it Степанец А.И.} Методы теории
приближений: В 2 ч. --- Киев: Ин-т математики НАН Украины, 2002.
--- Ч.\,1. --- 427\,c.

\item\label{StepanetsSerdyukShydlich2008} {\it Степанец А.И., Сердюк А.С., Шидлич А.Л.} Классификация бесконечно дифференцируемых периодических функций // Укр. мат. журн. - 2008. - 60, № 12. - С. 1686–1708.

\item\label{Sharapudinov1983} {\it Sharapudinov I.I.} On the best approximation and polynomials of the least quadratic deviation // Anal. Math. - 1983. - 9, № 3. - P.~223–234.

\item\label{Oskolkov1986} {\it Oskolkov K.I.} Inequalities of the "large sieve"\ type and applications to problems of trigonometric approximation // Anal. Math. - 1986. - 12, № 2. - P.~143–166. 

\item\label{Motornyj1990} {\it Моторний В.П.} Приближение периодических функций интерполяционными многочленами в $L_1$ // Укр. мат. журн. - 1990. - 42, № 6. - С. 781–786.
 
\item\label{Serdyuk2000} {\it Сердюк А.С.} Наближення періодичних функцій високої гладкості інтерполяційними тригонометричними поліномами в метриці $L_1$. - Укр. мат. журн. – 2000. – Т.52, № 7. – C.~994–998.
 
\item\label{Serdyuk2001} {\it Сердюк А.С.} Наближення інтерполяційними тригонометричними поліномами нескінченно диференційовних періодичних функцій в інтегральній метриці. - Укр. мат. журн. – 2001. – Т.53, № 12. – C.~1654–1663.
 
\item\label{Serdyuk2002} {\it Сердюк А.С.} Наближення періодичних аналітичних функцій інтерполяційними тригонометричними поліномами в метриці простору $L.$ - Укр. мат. журн. – 2002. – Т.54, № 5. – C.~692–699. 
 
\item\label{SerdyukVojtovych2010} {\it Сердюк А.С., Войтович В.А.} Наближення класів цілих функцій інтерполяційними аналогами сум Валле Пуссена // Збірник праць Інституту математики НАН України. - Т.7, № 1: Теорія наближення функцій та суміжні питання.- Київ: Ін-т математики НАН України, 2010.- C.~274-297.
 
\item \label{Kornejchuk1987} {\it Корнейчук Н.П.} Точные константы
в теории приближения. --- М.: Наука, 1987. --- 423~с.

\item\label{Serdyuk2005} {\it Сердюк А.С.} Наближення класів аналітичних функцій сумами Фур'є в рівномірній метриці //Укр. мат. журн. - 2005. - 57, № 8. - С.~1079 – 1096.

\item\label{Stechkin1980} {\it Стечкин С.Б.} Оценка остатка ряда Фурье для дифференцируемых функций // Приближение функций полиномами и сплайнами, Сборник статей, Тр. МИАН СССР. -  1980. - 145. - C.~126–151.

\item\label{Telyakovskii1989} {\it Теляковский С.А.} О приближении суммами Фурье функций высокой гладкости // Укр. матем. журн.~- 1989. - 41, № 4. - С.~510–518.

\end{enumerate}

\end{document}